\documentclass[12pt]{amsart}

\usepackage{amsmath}
\usepackage{amsfonts}
\usepackage{amssymb}
\usepackage{amsthm}
\usepackage{hyperref}
\usepackage{enumerate}
\usepackage{cleveref}
\usepackage{mathrsfs}
\usepackage{lscape}
\usepackage{geometry}
\usepackage{multicol}
\usepackage{todonotes}

\newtheorem{theorem}{Theorem}[section]
\newtheorem{lemma}[theorem]{Lemma}
\newtheorem{proposition}[theorem]{Proposition}

\Crefname{conjecture}{Conjecture}{Conjectures}

\theoremstyle{plain}

\theoremstyle{plain}

\newtheorem{remark}[theorem]{Remark}
\newtheorem*{remark*}{Remark}

\newcommand{\N}{\mathbb{N}}

\newcommand{\Z}{\mathbb{Z}}

\newcommand{\Q}{\mathbb{Q}}
\newcommand{\R}{\mathbb{R}}
\newcommand{\C}{\mathbb{C}}

\newcommand{{\D}}{\delta}

\newcommand{\SL}{\operatorname{SL}}

\renewcommand{\Im}{\operatorname{Im}}

\newcommand{\SLZ}{\SL_2(\Z)}

\newcommand{\calE}{\mathcal{E}}

\renewcommand{\d}{\partial}

\newcommand{\HH}{\mathfrak{H}}

\newcommand{\lcm}{\operatorname{lcm}}

\newcommand{\frakZ}{\mathfrak{Z}}

\newcommand{\zetahat}{\widehat\zeta}

\newcommand{\smallmat}[4]{\left(\begin{smallmatrix} #1 & #2 \\ #3 & #4\end{smallmatrix}\right)}
\newcommand{\fraka}{\mathfrak a}
\renewcommand{\d}{\mathrm d}
\newcommand{\frakz}{\mathfrak z}
\newcommand{\frakp}{\mathfrak p}

\numberwithin{equation}{section}
\numberwithin{table}{section}

\begin{document}
\author{Claudia Alfes and Michael H.\ Mertens}

\title{On Kleinian mock modular forms}

\thanks{The first author is partially supported by the Daimler and Benz Foundation, and funded by the Deutsche Forschungsgemeinschaft (DFG, German Research Foundation) -- SFB-TRR 358/1 2023 -- 491392403}

\address{Universität Bielefeld, Fakultät für Mathematik, Postfach 100 131, 33501 Bielefeld, Germany,
E-Mail: \url{alfes@math.uni-bielefeld.de}}

\address{Universit\"at zu K\"oln, Department Mathematik/Informatik, Abteilung Mathematik, Weyertal 86--90, 50931 K\"oln, Germany, E-Mail: \url{mmertens@math.uni-koeln.de}}

\begin{abstract}
We give an explicit and computationally efficient construction of harmonic weak Maass forms which map to weight $2$ newforms under the $\xi$-operator.  Our work uses a new non-analytic completion of the Kleinian $\zeta$-function from the theory of Abelian functions.
\end{abstract}

\maketitle

\section{Introduction and Statement of Results}
Harmonic weak Maass forms are real-analytic generalizations of classical modular forms which were introduced by Bruinier and Funke in \cite{brfu04}. By now harmonic weak Maass forms are ubiquitous in number theory and many other areas of mathematics and theoretical physics (see for instance \cite{onovisions,kenbook} and the references therein). A \textit{harmonic weak Maass form of weight $k\in\Z$ for a congruence subgroup $\Gamma_0(N)$} is a smooth function on the upper half-plane $\HH$ which transforms like a usual (holomorphic) modular form of weight $k$ under $\Gamma_0(N)$. Rather than being holomorphic, it is annihilated by the weight $k$ hyperbolic Laplace operator 
$$\Delta_k= -y^2\left(\frac{\partial^2}{\partial x^2}+\frac{\partial^2}{\partial y^2}\right)+iky\left(\frac{\partial}{\partial x}+i\frac{\partial}{\partial y}\right),$$
where we write $\tau=x+iy$ for $z\in\HH$.  In addition, they need to satisfy a certain growth condition at the cusps.

One of the central tools in the theory of harmonic weak Maass forms is the $\xi$-operator defined by $\xi_k:=- 2i y^k \overline{\frac{\partial}{\partial\overline{\tau}}}$.  As Bruinier and Funke first showed in \cite{brfu04}, this operator yields a surjective map from the space $H_k(N)$ of harmonic weak Maass forms of weight $k$ for $\Gamma_0(N)$ to the space $S_{2-k}(\Gamma_0(N))$ of (holomorphic) cusp forms of dual weight $2-k$. The image of a harmonic weak Maass form under the $\xi$-operator is called the \textit{shadow} of its canonical holomorphic part, which is itself referred to as a \emph{mock modular form}.

As there are infinitely many preimages of a given cusp form under the $\xi$-operator, we would like to identify \textit{distinguished} preimages. This can be achieved by employing Poincar\'e series \cite{Br02,BO07} or holomorphic projection \cite{ARZ,DMZ,MOR}. In recent work Ehlen, Li, and Schwagenscheidt were able to construct so-called \textit{good} preimages of CM forms using a certain theta lift. In particular, the holomorphic parts of such preimages have algebraic Fourier coefficients at $\infty$. Previous work of Bruinier, Ono, and Rhoades \cite{BOR} guarantees the existence of such good preimages.  Ehlen, Li, and Schwagenscheidt were able to determine the exact algebraic number field containing their Fourier coefficients.


For weight $2$ newforms associated to rational elliptic curves, i.e.\ with rational Fourier coefficients, the first author together with Griffin, Ono, and Rolen \cite{agor} constructed distinguished preimages under the $\xi$-operator extending earlier work of Guerzhoy \cite{guerzhoy}. This construction uses a lattice-invariant completion of the Weierstrass $\zeta$-function from the classical theory of elliptic functions. When evaluated at the Eichler integral of the newform this gives essentially a harmonic weak Maass form of weight $0$.

These results were recently extended to newforms with rational Fourier coefficients of positive weight by the authors in joint work with Funke and Rosu \cite{afmr}. 

In this paper, we extend the construction from \cite{agor} to newforms with non-rational coefficients. This directly leads to the question for an analogue of the Weierstrass $\zeta$-function in the context of Abelian functions, as we shall explain in the following paragraphs (see also Section~\ref{secklein} for further details).
 
Let $K/\Q$ be a number field and let $f\in S_2(N)$ be a newform of weight $2$ for $\Gamma_0(N)$ with coefficients in $K$ with Galois conjugates $f_1,=f,f_2,\ldots,f_r$. By $V$ we denote the associated component of the modular curve $X_0(N)$ and write $\Lambda_V=\omega\Z^r+\omega'\Z^r\subset \C^r$ for the associated period lattice. As it turns out (see Section~\ref{seccomputational}) we can choose $\omega,\omega'$ such that $\Omega=\omega^{-1}\omega'\in \HH_r$ is in the Siegel upper half-space of genus $r$ (see \eqref{eqsiegel}). 

For $u\in \C^r$ and $\alpha, \beta\in \R^r$ we define the Kleinian $\sigma$-function by
\[
\sigma\left[\begin{smallmatrix} \alpha \\ \beta\end{smallmatrix}\right](u;\Lambda_V)=\exp\left(\frac 12 u^{tr}\omega^{-1}\eta u\right)\theta\left[\begin{smallmatrix}\alpha\\ \beta\end{smallmatrix}\right]\left(\omega^{-1}u;\omega^{-1}\omega'\right).
\]
where $\theta\left[\begin{smallmatrix}\alpha\\ \beta\end{smallmatrix}\right]\left(u;\omega^{-1}\omega'\right)$ is the Riemann theta function of characteristic $\left[\begin{smallmatrix} \alpha \\ \beta\end{smallmatrix}\right]$ (compare \eqref{eq:theta}) and $\eta,\eta'$ denote the quasi-periods of $V$.  For brevity we sometimes write $\sigma(u)=\sigma\left[\begin{smallmatrix} \alpha \\ \beta\end{smallmatrix}\right](u;\Lambda_V)$.

We then define $\zeta(u):=\nabla_u \log\sigma(u)$. This function is analytic but not Abelian, i.e. invariant under translations by lattice points $\ell\in\Lambda_V$. Following an idea of Rolen \cite{Rolen15} we find a non-meromorphic completion $\widehat\zeta$ of the Kleinian $\zeta$-function in Proposition~\ref{propzeta}, which indeed satisfies
\[
\widehat\zeta(u+\ell)=\widehat{\zeta}(u) \text{ for all } \ell=\lambda\omega+\mu\omega'\in \Lambda_V
\]
wherever it is defined.

We define the vector of Eichler integrals associated to $\vec f =(f_1,\ldots, f_r)$ as
\[
\vec\calE(\tau) = \left( -2\pi i \int_\tau^{i\infty} f_1(z)dz,\ldots, -2\pi i \int_\tau^{i\infty} f_r(z)dz\right).
\]
The \textit{obstruction to modularity} of $\vec\calE(\tau) $ is an element of $\Lambda_V$, that is $\vec\calE(\gamma.\tau)  -\vec\calE(\tau)  \in \Lambda_V$, for $\gamma\in\Gamma_0(N)$. Evaluating the Kleinian $\zeta$-function at $\vec\calE(\tau) $ gives a distinguished preimage of $\vec f$.

\begin{theorem}\label{thmmain}
Let the notation be as above. The function 
$$\widehat \frakZ_V:\HH\to \C^{1\times r},\quad \tau\mapsto \widehat\zeta(\vec\calE(\tau))$$
is defined for all $\tau\in\HH$ such that the Kleinian $\sigma$-function $\sigma(\vec\calE(\tau);\Lambda_V)$ does not vanish. There it is $\Gamma_0(N)$-invariant and is annihilated by the hyperbolic Laplacian $\Delta_0$. Moreover, we have
$$\xi_0\widehat\frakZ_V(\tau)=4\pi^2 {\vec f}^{tr}P^{-1}$$
where
$$P=\frac{1}{2i}\left(\overline\omega\omega'^{tr}-\overline{\omega'}\omega^{tr}\right)$$
is positive definite.
\end{theorem}
For later reference we denote the meromorphic part of $\widehat \frakZ_V$ by
\begin{gather}\label{eqZV}
\frakZ_V(\tau)=\zeta(\vec\calE(\tau))-\frac 12 {\vec\calE(\tau)}^{tr}\left(\omega^{-1}\eta+\eta^{tr}\omega^{-tr}\right)+\pi \vec\calE(\tau)^{tr}P^{-1}
\end{gather}
so that
$$\widehat \frakZ_V(\tau)=\frakZ_V(\tau)-\pi \overline{\vec\calE(\tau)}^{tr}P^{-1}.$$ 
We call $\frakZ_V(\tau)$ the \emph{(polar) Kleinian mock modular form} associated to $V$.

\begin{remark}
For $r=1$, i.e.\ a newform with rational coefficients, \Cref{thmmain} yields
$$\xi_0 \widehat\frakZ_V(\tau)=\frac{4\pi^2}{\Im(\overline{\omega}\omega')}f(\tau)=\frac{4\pi^2}{\operatorname{vol}(\Lambda_V)}f(\tau),$$
recovering the corresponding result in \cite{agor}.
\end{remark}

The Kleinian mock modular form $\frakZ_V(\tau)$ may have poles. This phenomenon also occurs in the setup of \cite{agor}, where the authors show that there is a modular function which cancels the poles. We obtain the following result in this direction.


\begin{theorem}\label{thmmodfctn}
Assume the notation as above and choose the characteristic of the Riemann theta function as the Riemann characteristic of the base point $\infty$.  Assume in addition that $V$ is an $(n,s)$-curve, i.e.\ there is a model of the form
\begin{gather}\label{eqnscurve}
y^n=x^s+\sum_{\substack{i,j\geq 0 \\ in+j s< ns}} c_{ij}x^iy^j,\quad c_{ij}\in\Q\text{ and }s>n\text{ are coprime}
\end{gather}
and that for $\tau\in\HH$ the vector $\vec\calE(\tau)$ is not contained in the theta divisor except for points in $\Lambda_V$ (see \eqref{eqthetadiv}).

Then the principal part at $\infty$ of the scalar-valued function $\frakz_V(\tau):=\frakZ_V(\calE(\tau))\cdot (1,...,1)^{tr}$ has coefficients in $K$ and there exists a modular function $F_V$ for $\Gamma_0(N)$ with algebraic Fourier coefficients, such that $\frakz_V+F_V$ is a harmonic weak Maa\ss~form.
\end{theorem}

\begin{remark}
\begin{enumerate}
\item We note that every elliptic or hyperelliptic curve is in particular an $(n,s)$-curve. Ogg \cite{Ogg74} famously proved that the modular curve $X_0(N)$ is hyperelliptic if and only if 
$$N\in\{22,23,26,28,29,30,31,33,35,37,39,41,46,47,50,59,71\},$$
so that Theorem~\ref{thmmodfctn} applies in those cases. Furthermore it is a well-known fact that every curve of genus $2$ is hyperelliptic, wherefore our result also applies whenever the length $r$ of the Galois orbit of the considered newform is $\leq 2$.
\item Our proof of Theorem \ref{thmmodfctn} relies on the algebraicity of the coefficients of the expansion of the Kleinian $\sigma$-function around $0$ (\cite[Theorem 3]{Nakayashiki}). To the authors' knowledge this is only known for $(n,s)$-curves explaining this technical assumption on $V$. 
\item The additional assumption that $\vec\calE(\tau)$ lies in the theta divisor only if $\vec\calE(\tau)\in \Lambda_V$ is also a technical one. It allows us to only consider the expansion of $\sigma$ around $0$ rather than other points. The algebraic properties of expansions of $\sigma$ around other points do not seem to have been investigated in the literature.
\item So-called Millson theta lifts of the functions $\frakz_V+F_V$ fall into the framework of results of Bruinier and Ono \cite{brono}, i.e.\ it is possible to relate the algebraicity of the Fourier coefficients of these lifts to the vanishing of the twisted central $L$-derivatives of the associated newform of weight $2$ (compare \cite{agor} for the case of newforms with rational Fourier coefficients). 
\end{enumerate}
\end{remark}

The following theorem gives the expansion of the Kleinian mock modular form at other cusps than $\infty$.
\begin{theorem}\label{thmAL}
In the situation of \Cref{thmmain} let $Q$ be an exact divisor of $N$ (i.e. $\gcd(Q,N/Q))=1$) and $W_Q$ the corresponding Atkin-Lehner involution. Denote by $\lambda_Q\in\{\pm 1\}$ the eigenvalue of the newform $f$ under $W_Q$, i.e.\ $f|W_Q=\lambda_Qf$. Further let 
$$L_Q(f):=-2\pi i \int_{W_Q^{-1}.\infty}^\infty f(t)\d t.$$
Then we have
$$\widehat \frakZ_V|W_Q(\tau)=\widehat\zeta(\lambda_Q(\vec \calE(\tau)-L_Q(\vec f))).$$
\end{theorem}

The paper is organized as follows. In Section 2 we introduce the Riemann theta function, define the Kleinian $\sigma$ and $\zeta$-function, and construct the completed Kleinian $\zeta$-function. We prove Theorem \ref{thmmain}, \ref{thmmodfctn}, and \ref{thmAL} in Section 3. In Section 4 we present some computational examples.

\section*{Acknowledgments}
We thank Eberhard Freitag, Jens Funke, Eugenia Rosu, Fredrik Str\"omberg, Bernd Sturmfels, Don Zagier, and David Zureick-Brown for interesting discussions, and Annika Burmester and Paul Kiefer for comments on an earlier version of this manuscript. The first author thanks the Hausdorff Institute where this work was initiated.

\section{Construction of Kleinian Abelian functions}\label{secklein}
In this section we construct the analog of the lattice-invariant Weierstrass $\zeta$-function in the Abelian case. We employ the approach using Riemann theta functions to construct the $\sigma$-function.

\subsection{The Riemann $\theta$-function}
We review the definition and some basic properties of the Riemann $\theta$-function.

We denote by
\begin{gather}\label{eqsiegel}
\HH_g:=\{\Omega\in\C^{g\times g}\: :\: \Omega^{tr}=\Omega,\ \Im(\Omega)>0\}
\end{gather}
the \emph{Siegel upper half-space} of genus $g$. Here, we write $A>0$ to indicate that a real symmetric matrix $A$ is positive definite.  

Let $\alpha,\beta\in\R^g$. For $u\in\C^g$ and $\Omega\in\HH_g$ we define the \emph{Riemann theta function of characteristic $\left[\begin{smallmatrix} \alpha \\ \beta\end{smallmatrix}\right]$} by
\begin{equation}\label{eq:theta}
\theta\left[\begin{smallmatrix}
\alpha\\ \beta
\end{smallmatrix}\right]( u;\Omega):=\sum_{ m\in\Z^g} e\left(( m+\alpha)^{tr}(u+\beta)+\frac 12 (m+\alpha)^{tr}\Omega (m+\alpha)\right),
\end{equation}
where $e(x)=\exp(2\pi i x)$. We write $\theta(u;\Omega)=\theta[0](u;\Omega)$.

This function is well-known to satisfy the following transformation properties (see e.g.\ \cite[pp. 123, 194, 195]{Mumford} or \cite[Theta Transformation Formula 8.6.1]{birkenhakelange}) 
\begin{align}
\theta\left[\begin{smallmatrix}
\alpha\\ \beta
\end{smallmatrix}\right](u+\lambda\Omega+\mu;\Omega)&=e\left(-\frac 12 \lambda^{tr} \Omega\lambda- u^{tr}\lambda -\lambda^{tr}\beta+\mu^{tr}\alpha\right)\theta\left[\begin{smallmatrix}
\alpha\\ \beta
\end{smallmatrix}\right](u;\Omega),\quad \lambda,\mu\in\Z^g, \label{eqthetaell}\\
\theta(u;\Omega+S)&=\theta(u;\Omega),\quad \text{if } S\in\Z^{g\times g}_{sym}\text{ even},\\
\theta(\Omega^{-1}u;-\Omega^{-1})&=(\det -i\Omega)^{1/2}e\left(\frac 12 u^{tr}\Omega^{-1} u\right)\theta(u;\Omega)\label{eqthetamod2}.
\end{align}
These transformations imply that $\theta(u;\Omega)$ is a Siegel-Jacobi form of weight and index $1/2$ (see \cite{YY}).  

\subsection{Kleinian $\zeta$-functions and their completions}\label{secklein2}
Following the classical treatments by Klein and Baker \cite{Klein,Baker} and the more modern discussions, e.g.\ in \cite{BEL,EA12}, we construct Abelian functions via the Riemann theta function.

Let $V$ be a non-singular algebraic curve of genus $g$ with period matrices $\omega,\omega'\in \C^{g\times g}$, so that the Jacobian of $V$ is isomorphic to the quotient $\C^g/\Lambda_V$ where
$$\Lambda_V=\{\omega m+\omega' n\: :\: m,n\in\Z^g\}.$$
Here and throughout, we choose the period matrices $\omega,\omega'$ such that $\Omega=\omega^{-1}\omega'\in\HH_g$ lies in the Siegel upper half-space.
Further let $\eta,\eta'\in\C^{g\times g}$ denote the quasi-period matrices of $V$.  For an arbitrary characteristic $\left[\begin{smallmatrix} \alpha \\ \beta\end{smallmatrix}\right]\in\R^{2g}$ and $ u=(u_1,...,u_g)^{tr}\in \C^g$ we define the \emph{Kleinian $\sigma$-function} of characteristic $\left[\begin{smallmatrix} \alpha \\ \beta\end{smallmatrix}\right]$ by
\begin{equation}\label{eqsigma}
\sigma\left[\begin{smallmatrix} \alpha \\ \beta\end{smallmatrix}\right](u)=\sigma\left[\begin{smallmatrix} \alpha \\ \beta\end{smallmatrix}\right](u;\Lambda_V)=\exp\left(\frac 12 u^{tr}\omega^{-1}\eta u\right)\theta\left[\begin{smallmatrix}\alpha\\ \beta\end{smallmatrix}\right]\left(\omega^{-1}u;\omega^{-1}\omega'\right).
\end{equation} 

\begin{remark}
\begin{enumerate}
\item The definition of the Kleinian $\sigma$-function  above differs from that in \cite{BEL,EA12} by a constant depending on the curve $V$.
\item In \cite{BEL,EA12}, the characteristic is fixed as the Riemann characteristic of the base point $\infty$. In particular, this implies that $\alpha,\beta$ are half-integral. We usually choose the characteristic with half-integer entries, but not necessarily corresponding to the base point $\infty$.
\item If $V=E$ is an elliptic curve, then the Kleinian $\sigma$-function of characteristic $\left[\begin{smallmatrix} 1/2 \\ 1/2\end{smallmatrix}\right]$ coincides, up to a constant factor, with the classical Weierstrass $\sigma$-function.  
\end{enumerate}
\end{remark}

We require the following result on the zeros of the Kleinian $\sigma$-function (see e.g. \ \cite[Corollary II.3.6]{Mumford} and \cite[p.1661]{EA12}).
\begin{lemma}\label{lemthetadiv}
Choosing $\left[\begin{smallmatrix}\alpha\\\beta \end{smallmatrix}\right]\in \frac 12\Z^{2g}$ as the Riemann characteristic of the base point $\infty$, we have that $\sigma(u)=0$ if and only if
\begin{gather}\label{eqthetadiv}
u=\pm\left(\int_\infty^{P_1}du+...+\int_\infty^{P_{g-1}}du\right)+\ell,
\end{gather}
where $(P_1,...,P_{g-1})\in\operatorname{Sym}^{g-1}(V)$, where $\operatorname{Sym}^k(V)$ denotes the $k$th symmetric power of $V$, and $\ell\in\Lambda$.
\end{lemma}

From now on we suppress the characteristic from the notation when it can be chosen arbitrarily.

We define the \emph{$i$th Kleinian $\zeta$-function} by
\[
\zeta_i(u)=\partial_{u_i}\log\sigma(u)
\]
and  the \emph{Kleinian $\zeta$-function} by
\[
\zeta(u):=(\zeta_1(u),...,\zeta_g(u))=\nabla_u \log\sigma(u).
\]
The \emph{Kleinian $\wp$-functions} $\wp_{ij}(u):=-\partial_{u_j}\zeta_i(u)$ are then Abelian functions with respect to the lattice $\Lambda_V$, i.e.\ they satisfy the transformation law
$$\wp_{ij}(u+\ell)=\wp(u),\quad\ell\in\Lambda_V,$$
wherever they are defined.
 The $\zeta$-functions however are not Abelian but rather satisfy 
\begin{gather}
\zeta(u+\omega m+\omega' n)=\zeta(u)+(m^{tr}\eta +n^{tr}\eta').
\end{gather}
Again, for $V=E$ an elliptic curve this reduces to the well-known transformation law of the Weierstrass $\zeta$-function.

It is a classical fact going back to Eisenstein that the Weierstrass $\zeta$-function admits a non-analytic completion which is invariant under translation by lattice points \cite{eisenstein}. We adapt Rolen's proof of this property \cite{Rolen15} to the setting of Kleinian $\zeta$-functions.

The following Lemma follows by a straightforward computation.

\begin{lemma}\label{lemraising}
Let $F:\C^g\times \HH_g\to \C$ be a smooth function which satisfies the following functional equation
$$F(u+\Omega\lambda+\mu;\Omega)=e\left(-2m\lambda^{tr} u+f(\Omega,\lambda,\mu)\right) F(u;\Omega),\qquad \lambda,\mu\in\Z^g,$$
for some real number $m$ and some fixed function $f:\HH_g\times\R^g\times \R^g\to\C$.

Then the function defined by
$$\widetilde{F}(u;\Omega):=\frac{1}{F(u;\Omega)}\left(i\nabla_u-4\pi m \Im(u)^{tr}\Im(\Omega)^{-1}\right) F(u;\Omega)$$
satisfies
$$\widetilde F(u+\Omega\lambda+\mu;\Omega)=\widetilde F(u;\Omega),\quad \text{for all }\lambda,\mu\in\Z^g,$$
wherever $F(u,\Omega)\neq 0$.
\end{lemma}

\begin{remark}
The differential operator 
\begin{gather}\label{eqcalY}
\mathcal Y_+:=i\nabla_u-4\pi m \Im(u)^{tr}\Im(\Omega)^{-1}
\end{gather}
in \Cref{lemraising} looks similar to the raising operator
$$Y_+=i\partial_z-4\pi m\frac{\Im z}{\Im \tau}$$
in the theory of Jacobi forms (see e.g. \cite{BS98}). This operator maps Jacobi forms of weight $k$ and index $m$ to non-holomorphic Jacobi forms of weight $k+1$ and index $m$. An extension of this operator in the context of Siegel-Jacobi forms is given in \cite{YY}. However, the operator $\mathcal Y_+$ in \eqref{eqcalY} does not respect the action of the symplectic group in the same way as  \eqref{eqthetamod2}. For the purpose of this paper, this is not required and the operators in \cite{YY} are not suitable for our setup.
\end{remark}

We obtain the following result which is analogous to that in \cite{Rolen15}.
\begin{proposition}\label{propzeta}
The function
$$\widehat \zeta(u)=\zeta(u)-\frac 12 u^{tr}\left(\omega^{-1}\eta+\eta^{tr}\omega^{-tr}\right)+2\pi i\Im(\omega^{-1}u)^{tr}\Im(\omega^{-1}\omega')^{-1}\omega^{-1}$$
is a non-meromorphic Abelian function for the lattice $\Lambda_V$, i.e.\ for any $\ell=\lambda\omega+\mu\omega'\in \Lambda_V$ we have
$$\widehat\zeta(u+\ell)=\widehat{\zeta}(u)$$
wherever both sides are defined.
\end{proposition}
\begin{proof}
For simplicity we let
$$\vartheta(u):=\theta\left[\begin{smallmatrix}\alpha\\ \beta\end{smallmatrix}\right]\left(\omega^{-1}u;\omega^{-1}\omega'\right).$$
It follows immediately from \Cref{lemraising} that the function
\begin{gather}\label{eqdifftheta}
\frac1{\vartheta(u)}\nabla_u\vartheta(u)+2\pi i \Im(\omega^{-1}u)^{tr}\Im(\omega^{-1}\omega')^{-1}\omega^{-1}
\end{gather}
is a non-meromorphic Abelian function where it is defined. From \eqref{eqsigma} we obtain that
$$\zeta(u)=\nabla_u\left(\frac 12 u^{tr}\omega^{-1}\eta u +\log \vartheta(u)\right)=\frac 12 u^{tr}(\omega^{-1}\eta +\eta^{tr}\omega^{-tr})+\frac{1}{\vartheta(u)}\nabla_u(\vartheta(u)),$$
so the claim follows.
\end{proof}

\section{Proofs of the main results}
In this section we prove the main results of the paper.

\begin{proof}[Proof of \Cref{thmmain}]
We first note that $\zeta(u)$ is defined whenever $u$ does not lie in the divisor of the Kleinian $\sigma$-function.

We now prove the $\Gamma_0(N)$-invariance. 
 Let $\gamma\in\Gamma_0(N)$ and $\tau\in\HH$. We then have
$$\vec\calE(\gamma.\tau)=-2\pi i\int_{\gamma.\tau}^\infty \vec f(t) \d t=-2\pi i\left(\int_\tau^\infty -\int_{\tau}^{\gamma.\tau}\right)\vec f(t) \d t.$$
Since all components of $\vec f$ are holomorphic cusp forms, we have
$$\int_\tau^{\gamma.\tau} \vec f(t)\d t=\int_{\fraka}^{\gamma.\fraka} \vec f(t)\d t$$
for an arbitrary cusp $\fraka$ of $\Gamma_0(N)$ (see e.g.\ \cite[Proposition 10.5]{Stein}). Since a path from $\fraka$ to $\gamma.\fraka$ lies in the cuspidal homology of $X_0(N)$, we find by definition of the period lattice $\Lambda_V$ that $\ell=-2\pi i\int_\tau^{\gamma.\tau}\vec f(t)dt\in \Lambda_V$. The completed Kleinian $\zeta$-function $\widehat\zeta(u;\Lambda_V)$ is invariant under translations by points in $\Lambda_V$ by \Cref{propzeta}. Therefore, it follows that the function
$$\widehat\frakZ_V(\gamma.\tau)=\widehat\zeta(\vec\calE(\gamma.\tau))=\widehat\zeta(\vec\calE(\tau)-\ell)=\widehat\zeta(\vec\calE(\tau))=\widehat\frakZ_V(\tau)$$
is indeed $\Gamma_0(N)$-invariant wherever it is defined.

We now proceed to show the properties under the action of the differential operators $\Delta_0$ and $\xi_0$.
Using \Cref{propzeta} and \eqref{eqdifftheta} we can write
\begin{align}
\widehat{\zeta}(u)&=\frac1{\theta(\omega^{-1}; \Omega)}\left(\nabla_u\theta\right)(\omega^{-1}u; \Omega)\omega^{-1}+2\pi i \Im(\omega^{-1}u)^{tr}\Im(\Omega)^{-1}\omega^{-1}\nonumber \\
&=\frac1{\theta(\omega^{-1}u;\Omega)}\left(\nabla_u\theta\right)(\omega^{-1}; \Omega)\omega^{-1}+\pi  u^{tr}\omega^{-tr}\Im(\Omega)^{-1}\omega^{-1}\nonumber \\
&\qquad\qquad\qquad\qquad\qquad\qquad\qquad\qquad\quad-\pi \overline{u}^{tr}\overline\omega^{-tr}\Im(\Omega)^{-1}\omega^{-1}.\label{eqzetahatdiff}
\end{align}
This immediately implies that
\[
\xi_0\left(\widehat\frakZ_V(\tau)\right)=4\pi^2\vec f^{tr}\omega^{-tr}\Im(\Omega)^{-1}\omega^{-1}.
\]
A straightforward computation gives
$$\omega\Im(\Omega)\omega^{tr}=\frac{1}{2i}(\overline \omega\omega'^{tr}-\overline\omega'\omega^{tr})=P.$$
Since $\Omega=\omega^{-1}\omega'\in\HH_r$ we have that
$$2\Im(\Omega)=i(\overline \Omega^{tr}-\Omega)=i(\overline{\omega'}^{tr}\overline{\omega}^{-tr}-\omega^{-1}\omega')$$
is positive definite. Therefore the same is true for
$$i\omega(\overline{\omega'}^{tr}\overline{\omega}^{-tr}-\omega^{-1}\omega')\overline\omega^{tr}=2\overline P$$
and hence $P$ is positive definite.

\end{proof}

The proof of \Cref{thmmodfctn} is in part analogous to that of the corresponding result in \cite{agor}. Since the proof given there is rather short, we give a more detailed version here.

\begin{proof}[Proof of \Cref{thmmodfctn}]
By our assumption on $\vec\calE$ and Lemma \Cref{lemthetadiv} we see that by construction $\frakZ_V$ has a pole in $\tau$ if and only if $\vec\calE(\tau)\in\Lambda_V$.  Since $\zetahat(u)$ is lattice invariant, it is therefore enough to consider the expansion of $\zeta(u)$ around $u=0$. By \cite[Theorem 3]{Nakayashiki} we have
$$\sigma(u)=S_{\lambda(n,s)}(u)+\text{higher order terms},$$
where $S_{\lambda(n,s)}$ is the so-called \emph{Schur function} associated to the curve $V$ (for a precise definition see p.\ 192 of loc.cit.). This polynomial has rational coefficients. Since each newform $f_1,...,f_r$ has coefficients in $K$, so do the functions $\calE_1,...,\calE_r$. Therefore, we see that upon plugging $\vec\calE(\tau)$ into $\partial_{u_j}\sigma(u)/\sigma(u)$, the principal part of $\zeta(\vec\calE(\tau))$ at $\infty$ has coefficients in $K$ as well.

Next we show the existence of the meromorphic modular function $F_V$ which cancels all the poles of $\frakz_V$ within the upper half-plane: It is well-known (see e.g. Section~\ref{secklein2}) that any partial derivative of the Kleinian $\zeta$-function yields (up to sign) a Kleinian $\wp$-function, thus a meromorphic Abelian function. Therefore the function
$$\frakp_V(\tau)=\sum_{j=1}^r\partial_{u_j}^2\log \sigma (\calE(\tau))$$ 
is a meromorphic modular function with respect to $\Gamma_0(N)$ by the same argument employed in the proof of \Cref{thmmain}. The poles of this function within the upper half-plane are clearly at the same points as those of $\frakz_V$, but strictly with higher order.  As indicated in \cite{agor} we follow the proof of \cite[Theorem 11.9]{Cox}, which states that every modular function for $\Gamma_0(N)$ is a rational function in $j(\tau)$ and $j(N\tau)$. Let $\gamma_1,...,\gamma_{\iota(N)}$, $\iota(N)=[\Gamma_0(1):\Gamma_0(N)]$, be a fixed set of representatives of $\Gamma_0(N)\backslash\SLZ$ and assume $\gamma_1=\smallmat 1001$. We consider the function
$$G(X,\tau)=\sum_{i=1}^{\iota(N)} \frakp_V(\gamma_i\tau)\prod_{j\neq i}(X-j(\gamma_j\tau)).$$
This is clearly a polynomial in $X$ whose coefficients are meromorphic functions in $\tau$. In fact it is not hard to show that these coefficients are modular functions for $\SLZ$, whence they are all rational functions in $j(\tau)$.  We may therefore write
\[
G(X,j(\tau))=\sum_{k=0}^{\iota(N)-1} \frac{p_k\left(j(\tau)\right)}{q_k\left(j(\tau)\right)}X^k
\]
for certain polynomials $p_k,q_k\in\C[Y]$. 

In fact we can choose $p_k,q_k$ with algebraic coefficients. By assumption $V$ is an $(n,s)$-curve, so it follows from \cite[Theorem 3]{Nakayashiki} that the coefficients of the Taylor expansion of $\sigma$, and therefore of the Laurent expansions of both $\frakz_V$ and $\frakp_V$ are rational polynomials in the curve coefficients $c_{ij}$ in \eqref{eqnscurve}, and hence algebraic. Since the newform $f$ has algebraic Fourier coefficients at all cusps, it also follows that the modular function $\frakp_V$ has algebraic Fourier coefficients at all cusps. 

Let 
$$Q=\lcm(q_1,...,q_{\iota(N)})=\prod_{\ell=1}^M(Y-\alpha_\ell)$$ 
for some $\alpha_\ell\in\overline \Q$ and $M\in\N_0$. Now arguing exactly as in the aforementioned proof of \cite[Theorem 11.9]{Cox}, we find that we can write
\begin{gather}\label{eqfrakpj}
\frakp_V(\tau)=\frac{\sum_{k=0}^{\iota(M)-1}\tilde p_k(j(\tau))j(N\tau)}{\prod_{\ell=1}^M(j(\tau)-\alpha_\ell)\cdot \prod_{m\neq 1}(j(N\tau)-j(N\gamma_m\tau))}.
\end{gather}
Note that the numerator in \eqref{eqfrakpj} is holomorphic in $\HH$ and each factor in the denominator yields a simple pole of $\frakp_V$ in $\HH$ (we ignore the slight technical complication of elliptic fixed points for the sake of simplicity).  By multiplying through by all but one of the factors in the denominator (after canceling against potential zeros in the numerator), we obtain a modular function with algebraic Fourier coefficients with a simple pole precisely where $\frakz_V$ has a pole. Thus, we can cancel all the poles using only modular functions with algebraic coefficients. 
\end{proof}

As the proof of \Cref{thmAL} is almost literally the same as that of the analogous result in \cite{agor} (Theorem 1.2) we omit it here.


\section{Examples}

\subsection{Computational aspects}\label{seccomputational}

We briefly outline how to compute the quantities required for the construction of the Kleinian mock modular forms. 

Most of the facts in this section are by now fairly standard and more or less implemented in computer algebra systems like Sage \cite{sage}, Magma \cite{magma}, or Pari/Gp \cite{pari}.  We loosely follow the accounts in \cite{Chai,Wang} and Kapitel VI of \cite{Freitag}. 

Let $f\in S_2(N)$ be a newform whose coefficients lie in a number field $K/\Q$ and let $f_1=f,f_2,...,f_r$ denote its Galois conjugates. The vector of all these conjugates is denoted by $\vec f=(f_1,...f_r)^{tr}$. Suppose we have Fourier expansions 
$$f_j(\tau)=\sum_{n=1}^\infty a_j(n)q^n,\quad q=e^{2\pi i\tau},\quad a_j(n)\in K.$$
Then there is a component over $\Q$ of the modular curve $X_0(N)$ associated to the Galois orbit of $f$. Its Jacobian is given by $\C^r/\Lambda_f$ for the period lattice $\Lambda_f.$ We can find a basis for this lattice by computing the integrals
$$-2\pi i \int_\gamma \vec{f}(z) \d z,$$
where $\gamma$ runs through a basis of the integral homology $H^1(X_0(N),\Z)$, which can in turn be determined using the available functions in Sage or Magma.  

This may be achieved very efficiently by evaluating holomorphic Eichler integrals
$$\calE_j(\tau):=-2\pi i\int_\tau^\infty f_j(z)\d z =\sum_{n=1}^\infty \frac{a_j(n)}{n}q^n$$
at suitable points $\tau$ in the upper half-plane. 

It follows from work of Hida \cite{Hida} together with standard linear algebra that we can choose a basis $(a_1,...,a_r,b_1,...,b_r)$ of $H^1(X_0(N),\Z)$ with the property that the cycles follow the intersection pattern
$$a_i\circ a_j =0,\quad b_i\circ b_j=0,\quad a_i\circ b_j=e_i\delta_{ij},$$
where $\delta_{ij}$ denotes the usual Kronecker delta and $e_1\mid e_2\mid ...\mid e_r$ are positive integers. With respect to this basis we obtain matrices 
$$\omega=-2\pi i\left(\int_{a_i}f_j(z) \d z\right)_{i,j=1,...r},\quad \omega'=-2\pi i\left(\int_{b_i}f_j(z)\d z\right)_{i,j=1,...,r}$$
with the property that $\Omega:=\omega^{-1}\omega'\in\HH_r$ lies in the Siegel upper half-space.  An algorithm to compute this basis was found by Merel \cite{Merel} and is implemented e.g.\ in Magma.

Note that since we have 
\[
\zetahat(u)=\frac{1}{\vartheta(u)}\nabla_u\vartheta(u)+\pi u^{tr}P^{-1}-\pi \overline{u}^{tr}P^{-1}
\]
with $\vartheta(u)=\theta\left[\begin{smallmatrix} \alpha\\ \beta\end{smallmatrix}\right](\omega^{-1}u,\omega^{-1}\omega')$ and $P$ as in \Cref{thmmain} by \eqref{eqzetahatdiff}, we do not require the quasi-periods $\eta,\eta'$ to compute the Kleinian mock modular form.


\subsection{Level 27}\label{sec27}
We consider the unique newform $f=\eta(3\tau)^2\eta(9\tau)^2\in S_2(27)$ associated to the elliptic curve
\[
y^2+y=x^3-7  \quad \text{(LMFDB label 27.a3)}.
\]
It has rational coefficients and complex multiplication by $\Q(\sqrt{-3})$.



Since $f$ has rational coefficients, the results in \cite{agor} apply and we can compute a mock modular form whose shadow is $f$ (up to a constant multiple).  Alternatively, we can apply the strategy of this paper and find that the period lattice of $f$ is generated by
$$\omega=-0.294439 - 0.509984i\quad\text{and}\quad\omega'=-1.01996i.$$
Using the Kleinian zeta function with characteristic $\alpha=\beta=1/2$ we we employ \Cref{thmmain} to construct the function
$$\frakZ_V(\tau)=q^{-1}+\frac 12 q^2-\frac{701}{5}q^5+\frac{1407}{4}q^8-\frac{40776}{11}q^{11}+\frac{37961}{2}q^{14}-\frac{2125098}{17}q^{17}+O(q^{20}),$$
whose shadow is $4\pi^2 f$.  Note that the Fourier coefficients above are indeed rational numbers, which can be shown using work of Bruinier-Ono-Rhoades \cite[Theorem 1.3]{BOR} or Ehlen-Li-Schwagenscheidt \cite[Corollary 1.2]{ELS}.

\subsection{Level 23}
The modular curve $X_0(23)$ has genus $2$ and there is one Galois orbit of newforms, generated by the form with Fourier expansion
$$f(\tau)=q-\phi q^2+(2\phi-1)q^3+(\phi-1)q^4-2\phi q^5+O(q^6),\quad \phi=\frac{1-\sqrt{5}}{2}.$$
We denote the Galois conjugate of $f$ by $f^\sigma$. 

The four elements $\{-1/19,0\},\{-1/17,0\},\{-1/15,0\},\{-1/11,0\}$ form a basis of $H_1(X_0(23),\Z)$. Consequently we find the following basis for the period lattice
\begin{align*}
c_1&=\begin{pmatrix}
-1.062972 + 2.060558i \\
 0.642714 + 0.710672i
\end{pmatrix},\ c_2=\begin{pmatrix}
1.062972 + 2.060558i \\
 -0.642714 + 0.710672i
\end{pmatrix},\\
c_3&=\begin{pmatrix}
1.719925 + 0.787063i\\
 0.397219 + 1.860563i
\end{pmatrix},  \ c_4=\begin{pmatrix}
1.313906\\ 
2.079867
\end{pmatrix}.
\end{align*}
Computing the intersection pairing with the help of Magma we compute the period matrices 
\begin{align*}
\omega&=(c_1-c_2+c_3,c_2)=\begin{pmatrix}
-0.406019 + 0.787063i & 1.062972 + 2.060558i\\
1.682647 + 1.860563i & -0.642714 + 0.710672i
\end{pmatrix},\\
 \omega'&=(-c_2+c_3-c_4,c_2-c_3)=\begin{pmatrix}
   -0.656953 - 1.273495i & -0.656953 + 1.273495i \\
-1.039933 + 1.149891i & -1.039933 - 1.149891i
\end{pmatrix},
\end{align*}
so that we have
$$\Omega=\omega^{-1}\omega'=\begin{pmatrix}
0.01074169 & -0.3817894 \\ -0.3817894 & 0.3888885
\end{pmatrix}+i\begin{pmatrix}
0.7666448 &-0.1730782 \\ -0.1730782& 0.6607763
\end{pmatrix} \in\HH_2.$$

Choosing the characteristic $\alpha=(1/2,0)^{tr}$, $\beta=(1/2,1/2)^{tr}$, we find that $\theta(0,\Omega)=0$, so we directly obtain from \Cref{thmmain} that the components of the function 
$$\frac{1}{4\pi^2}\frakZ_{X_0(23)}\omega\Im(\Omega)\overline\omega^{tr}$$
yield preimages of the newforms $f$ and $f^\sigma$ resp.\ under $\xi_0$, up to the addition of a meromorphic modular form. Their meromorphic parts are given by
\begin{multline*}
0.259008q^{-1} + 1.000942 + 4.868978q + 18.294037q^2 + 68.247223q^3 + 252.912538q^4 \\
+ 938.377980q^5 + 3477.898343q^6 + 12892.503560q^7 + 47787.961740q^8+O(q^9)
 \end{multline*}
 and 
 \begin{gather*}
 -0.505669q^{-1} - 1.954167 - 6.9786217q - 26.191387q^2 - 97.573609q^3\\
  - 361.535343q^4
  - 1341.254086q^5 - 4971.053026q^6 - 18427.581035q^7\\
   - 68304.578170q^8+O(q^9).
 \end{gather*}
 
We consider the sum of the components of the vector-valued function. This will be a scalar-valued polar mock modular form $\frakz_V(\tau)$ whose shadow is some linear combination of the newforms $f$ and $f^{\sigma}$.  Note that the coefficient of $q^{43}$ and that of $q^{109}$ in $f$ both vanish, so it follows from the work of Bruinier-Ono-Rhoades \cite{BOR}, that the coefficients of $q^{43}$ and $q^{109}$ of a \emph{good} preimage of $f$ under $\xi_0$ should be algebraic numbers. Even though $\frakz_V$ is not guaranteed to be a good preimage in the sense of \cite{BOR}, we still find that 
 \begin{multline*}
 \frakz_V(\tau)=q^{-1}+3.864515+0.142266q+0.319448q^2+0.193313q^3+0.304709q^4\\
 +0.055558q^5+0.059060q^6+0.080332q^7+0.572492q^8-0.190607 q^9+O(q^{10})
 \end{multline*}
 and the coefficient of $q^{43}$ is (within computational precision) $27/43$ and that of $q^{109}$ is $942/109$.

We conclude this example by mentioning a few numerical observations. 
\begin{enumerate}
\item In this particular case, we see that the matrix $P=\frac{1}{2i}\left(\overline\omega\omega'^{tr}-\overline{\omega'}\omega^{tr}\right)$ from \Cref{thmmain} is diagonal, in fact we have, up to computational precision
$$P=\begin{pmatrix}
 3.741508 & 0 \\
 0 & 5.347829
 \end{pmatrix}=4\pi^2\begin{pmatrix}
 \langle f,f\rangle & 0 \\
 0 & \langle f^\sigma,f^\sigma\rangle
 \end{pmatrix},$$
 where $\langle\cdot,\cdot\rangle$ denotes the Petersson inner product. Possibly, this is a consequence of Haberland's formula for subgroups (see e.g.\ \cite[Theorem 5.2]{Cohen}).
\item
By the Petersson coefficient formula we can write the (conditionally convergent) cuspidal Poincar\'e series $\mathcal P(2,1,23;\tau)=\sum_{\gamma\in \Gamma_\infty\backslash\Gamma_0(23)} e^{2\pi i\tau}|_2 \gamma$ as
$$\mathcal P(2,1,23;\tau)=\frac{1}{4\pi\langle f,f\rangle} f(\tau)+\frac{1}{4\pi\langle f^{\sigma},f^{\sigma}\rangle} f^\sigma(\tau).$$
It is well-known that the preimage of a cuspidal Poincar\'e series under the $\xi$-operator is given by the so-called Maass-Poincar\'e series of dual weight, denoted by $\mathcal{Q}(2,-1,23;\tau)$ (see e.g.\ \cite[Theorem 6.10]{kenbook}). Computing the Fourier expansion of this Poincar\'e series numerically (see \cite[Theorem 6.10]{kenbook} for a description of the coefficients) strongly suggests that indeed
$$\widehat{\frakz}_V(\tau)=\mathcal{Q}(2,-1,23;\tau)+C$$
for some constant $C$, which would imply that $\widehat{\frakz}_V$ indeed has no poles within the upper half-plane.  Since their shadows are equal, we know that the difference  $\widehat{\frakz}_V(\tau)-\mathcal{Q}(2,-1,23;\tau)$ is a meromorphic modular function. By analyzing the behavior at cusps (see \Cref{thmAL}), we see that this function must have all its poles in the upper half-plane.

\item
Since $f|W_{23}=-f$, it follows by construction that $\frakZ_V(\tau)+\frakZ_V|W_{23}(\tau)$ should be a meromorphic modular function for $\Gamma_0(23)$ or rather the group $\Gamma_0(23)^+$. Indeed we find,  within computational accuracy, that
\begin{align*}
 &\frakZ_V(\tau)+\frakZ_V|W_{23}(\tau)\\
 \quad &=\begin{pmatrix} C(q^{-1}+\alpha+4q+7q^2+13q^3+19q^4+33q^5+47q^6+O(q^7))\\
C'(q^{-1}+\alpha'+4q+7q^2+13q^3+19q^4+33q^5+47q^6+O(q^7))
\end{pmatrix}
\end{align*}
for constants $C=-2.732921...,C'=3.7329211...,\alpha=- 0.019847,\alpha'= 2.543165...$. Note that the coefficients given above agree, apart from the constant term, with those of the Hauptmodul for the group $\Gamma_0(23)^+$ given by
$$T_{23+}(\tau)=t(\tau)+4\frac{t(\tau)}{(t(\tau)-1)},\quad t(\tau)=\frac{\eta(\tau)\eta(23\tau)}{\eta(2\tau)\eta(46\tau)}$$
(see e.g.\ \cite[Tables 3 and 4a]{CN79}, correcting an error in loc.cit.).
\end{enumerate}

\subsection{Level 256}
We consider the newform 
\[
f(\tau)=q+2\sqrt 2 q+5q^9-2\sqrt 2q^{11}+6q^{17}-6\sqrt{2}q^{19}+O(q^{21})\in S_2(256),
\]
which has CM by $\Q(\sqrt{-2})$. This is the smallest level for which there exists a CM newform with non-rational coefficients. As in the case of level $27$ in \Cref{sec27} the coefficients of the Kleinian mock modular form are algebraic.  As before we denote the Galois conjugate of $f$ by $f^\sigma$.



In this case one may check (e.g.\ by going through a list of generators of $\Gamma_0(256)$) that the associated period lattice in $\C^2$ is generated by
\begin{align*}
c_1&=-2\pi i\int_{49/512}^\infty \vec f(t)\d t=\begin{pmatrix} 2.767505- 2.767505i \\ 1.146338+ 1.146338i\end{pmatrix},\\
c_2&=-2\pi i\int_{55/512}^\infty \vec f(t)\d t=\begin{pmatrix}
1.956921- 1.956921i \\ 
-0.810583- 0.810583i
\end{pmatrix},\\
c_3&=-2\pi i\int_{29/256}^\infty \vec f(t)\d t=\begin{pmatrix}
0.810583+ 1.956921i\\ 1.956921+ 0.810583i
\end{pmatrix},\\
c_4&=-2\pi i\int_{39/256}^\infty \vec f(t)\d t=\begin{pmatrix}
1.956921+ 1.956921i \\
-0.810583+ 0.810583i
\end{pmatrix}.
\end{align*}
A suitable basis as described in \Cref{seccomputational} is given as follows
\begin{align*}
\omega&=(c_1,c_2-c_4)=\begin{pmatrix}
2.767505- 2.767505i &  - 3.913843i \\
1.146338+ 1.146338i &  -1.621166i
\end{pmatrix},\\
\omega'&=(c_2+c_3+2c_4,c_2+c_3+c_4)\\
&=\begin{pmatrix}
6.681348+ 3.913843i &  4.724426+ 1.956921i\\
-0.474828+ 1.621166i &  0.335754+ 0.810583i
\end{pmatrix},
\end{align*}
yielding
$$\Omega=\omega^{-1}\omega'=i\begin{pmatrix} 
1.414213 &  0.707106 \\
0.707106 & 0.707106
\end{pmatrix}=i\frac{\sqrt{2}}{2}\begin{pmatrix}
2 & 1 \\ 1 & 1
\end{pmatrix}\in\HH_2.$$
The latter equality is to be understood within computational precision (100 digits),  the authors are not aware of a rigorous proof for this observation.

The sum of the components $\frakz_V$ of $\frakZ_V$ is then given by
\begin{align*}
\frakz_V(\tau)&=q^{-1}+(1.828427+ 0.585786i)q + (-0.828427 - 3.313708i)q^3\\
&\qquad\qquad\qquad + (-5.254833 - 7.484271i)q^5 + (-7.487074 + 45.254833i)q^7 \\
&\qquad\qquad\qquad\qquad\qquad + (27.456710 + 44.496608i)q^9+O(q^{11})\\
&=q^{-1}+[(-1+2\sqrt 2)+i(2-\sqrt 2)]q +[(2-2\sqrt2)+i(8-8\sqrt 2)]q^3\\
&\quad +\frac 15[(200-160\sqrt 2)+i(104+100\sqrt 2)]q^5 +\frac 17[(1441-1056\sqrt 2)+i224\sqrt 2]q^7\\
&\quad\quad +\frac 19[(5211-3510\sqrt 2)+i(-5238+3987\sqrt 2)]q^9+O(q^{11}).
\end{align*}
Again, the last equality is to be understood within computational accuracy.

\bibliographystyle{plain}
\bibliography{bib.bib}

\end{document}